\documentclass[12pt]{amsart}
\usepackage[english]{babel}
\usepackage{amssymb,amsmath,amscd,amsthm}

\usepackage{cite}
\def\?{}

\newtheoremstyle{neosn}{0.5\topsep}{0.5\topsep}{\rm}{}{\sc}{.}{ }{\thmname{#1}\thmnumber{ #2}\thmnote{ {\mdseries#3}}}
\theoremstyle{neosn}

\newcommand{\GF}{\,\mathrm{GF}\,}
\newcommand{\ZZ}{{\mathbb Z}}

\begin{document}

\begin{center}
\large{\textbf{Centrally Essential Group Algebras}}
\end{center}

\hfill {\sf V.T. Markov}

\hfill Lomonosov Moscow State University

\hfill e-mail: vtmarkov@yandex.ru

\hfill {\sf A.A. Tuganbaev}

\hfill National Research University "MPEI"

\hfill Lomonosov Moscow State University

\hfill e-mail: tuganbaev@gmail.com

{\bf Abstract.} A ring $R$ with center $C$ is said to be \textit{centrally essential} if the module $R_C$ is an essential extension of the module $C_C$. In the paper, we study groups whose group algebras over fields are centrally essential rings. We focus on the centrally essential modular group algebras of finite groups over
fields of nonzero characteristic.

V.T.Markov is supported by the Russian Foundation for Basic
Research, project 17-01-00895-A. A.A.Tuganbaev is supported by
Russian Scientific Foundation, project 16-11-10013.


\begin{center}
\textbf{1. Introduction}
\end{center}
All considered rings are associative unital rings. A ring $R$ with
center $C$ is said to be \emph{centrally essential} if the module
$R_C$ is an essential extension of the module $C_C$. Any
commutative ring is a trivial example of a centrally essential
ring. However, there also exist noncommutative centrally essential
rings; in particular, there are noncommutative centrally essential
group rings of finite groups. For example, let $F=GF(2)$ be the
field of order 2 and $G=Q_8$ the quaternion group, i.e., $G$ is
the group with two generators $a,b$ and defining relations
$a^4=1$, $a^2=b^2$ and $aba^{-1}=b^{-1}$; see \cite[Section
4.4]{Hall}. Then the group algebra $FG$ is a noncommutative
centrally essential ring consisting of
256 elements (this follows from Theorem 1.1 below).

The main result of the paper is Theorem 1.1.

\textbf{Theorem 1.1.} Let $F$ be a field of characteristic $p>0$.

\textbf{1.} If $G$ is an arbitrary finite group, then the group algebra $FG$ is a centrally essential ring if and only if $G=P\times H$, where $P$ is the unique Sylow $p$-subgroup of the group $G$, the group $H$ is commutative, and the ring $FP$ is centrally essential.

\textbf{2.} If $G$ is a finite $p$-group and the nilpotence class\footnote{It is well known that every finite $p$-group is nilpotent, e.g., see \cite[Theorem 10.3.4]{Hall}.} of $G$ does not exceed $2$, then the group algebra $FG$ is a centrally essential ring.

\textbf{3.} There exists a group $G$ of order $p^5$ such that the group algebra $FG$ is not a centrally essential ring.

\textbf{Remark 1.2.} In connection to Theorem 1.1, we note that for an arbitrary field $F$ of zero characteristic and every group $G$, the group algebra $FG$ is a centrally essential ring if and only if the algebra $FG$ is commutative; see Proposition 3.2 in Section 3. Therefore, when studying centrally essential group algebras over fields, only the case of fields of positive characteristic is interesting.

The proof of Theorem 1.1 is given in the next section. We give
some necessary notions.

For an arbitrary finite subset $S$ of the group $G$ and any ring
$A$, we denote by  $\Sigma_S$ the element $\sum_{x\in S}x$ of the
group ring $AG$.

We use the following notation. Let $R$ be a ring and $G$ a group.
The center of the ring $R$ is denoted by $C(R)$. For any two
elements $a,b$ of the ring $R$, we set $[a,b]=ab-ba$. Following
\cite{Hall}, we set $(x,y)=x^{-1}y^{-1}xy$  for any two elements
$x,y$ of the group $G$; additive commutators and multiplicative
commutators are designated differently, since the elements of the
group are also considered as elements of the group ring. For any
element $g$ of the group $G$, we denote by $g^G$ the class of
conjugate elements which contains $g$. For a group $G$, the
\textit{upper central series} of $G$ is the chain of subgroups
$\{e\}=Z_0(G)\subseteq Z_1(G)\subseteq\ldots$,  where
$Z_i(G)/Z_{i-1}(G)$ is the center of the group $G/Z_{i-1}(G)$,
$i\geq 1$. We denote by $NC(G)$ the \textit{nilpotence class} of
the group $G$, i.e., the least positive integer $n$ with
$Z_n(G)=G$ (if it exists).

\begin{center}
\textbf{2. The proof of Theorem 1.1}
\end{center}

\textbf{Lemma 2.1.}\label{gab}
Let $A$ be a ring and $G$ a group. If the group ring $R=AG$ is
centrally essential, then $A$ also is a centrally essential ring
and the group $G$ is an $FC$-group, i.e., all the classes of
conjugate elements in $G$ are finite.

\textbf{Proof.}
Let $0\neq a\in A$ and $0\neq ca\in C(R)$. We have $c=\sum_{g\in
G}c_gg$. It follows from the relations $0=[c,b]=\sum_{g\in
G}[c_g,b]g$ for any $b\in A$ that $c_g \in C(A)$ for all $g\in G$.
Similarly, we have $ac_g\in C(A)$ for any $g\in G$. Since there
exists at least one element $g\in G$ with $ac_g\neq 0$, we obtain
our assertion about the ring $A$.

Now let $g\in G$. It is well known (e.g., see \cite[Lemma
4.1.1]{Passman}) that $C(AG)$ is a free $C(A)$-module with basis
$$
\{\Sigma_K\;|\;K\mbox{\ is a finite the class of of conjugate elements in the group\ }G\}.
$$
Therefore, if $0\neq cg=d$, where $c,d\in  Z(AG)$, then we can
compare the coefficients of $g$ in the left part and the right
part of the relation $cg=d$ and obtain that there exist finite
classes of conjugate elements $K_1,K_2$ such that $gx=y$ for some
$x\in K_1, y\in K_2$. For any $h\in G$, we have that
$hgh^{-1}=(hxh^{-1})^{-1}hyh^{-1}$ and the number of the right
parts of this relation does not exceed $|K_1|\cdot|K_2|$;
consequently, the class of elements, which are conjugate to $g$,
is finite.~\hfill$\square$

\textbf{Lemma 2.2.}\label{gab}
Let $A$ be a ring and $G$ a commutative monoid. If $A$ is a
centrally essential ring, then the monoid ring $AG$ is centrally
essential.

\textbf{Proof.}
Let $0\neq r=\sum_{g\in G}r_gg\in R$. We use the induction on the
number $k=k(r)$ of the coefficients $r_g$ which are not contained
in the center of the ring $A$. If $k(r)=0$, i.e., $r_g\in C(A)$
for all $g\in G$, then it is nothing to prove. Otherwise, we take
an element $h\in G$ with $r_h\not\in C(A)$.  There exists an
element $c\in C(A)$ with $0\neq cr_h\in C(A)$. Then it is clear
that $0\neq cr=\sum_{g\in G}cr_gg$ and $k(cr)<k(r)$. By the
induction hypothesis, there exists an element $d\in C(R)$ with
$0\neq dcr\in C(R)$. Since $dc\in C(R)$, the proof is
finished.~\hfill$\square$

\textbf{Lemma 2.3}\cite{MT}. In any centrally essential ring $R$, the following quasi-identity system
\begin{multline}\label{eqquid}
\forall n\in\mathbb{N}, x_1,\ldots,
x_n,y_1,\ldots,y_n,r\in R,\\
\left\{\begin{array}{rcl}x_1y_1+\ldots+x_ny_n&=&1\\
x_1ry_1+\ldots+x_nry_n&=&0\end{array}\right.
\Rightarrow r=0.
\end{multline}
is true. In particular, all idempotents of the centrally essential
ring $R$ are central.

\textbf{Proof.} We assume that $R$ is a centrally essential ring such that the relations from the above quasi-identity system hold, but $r\ne 0$. Then there exist two elements $c,d\in C(R)$ with
$cr=d\neq 0$. Consequently,
$$d=d(x_1y_1+\ldots+x_ny_n)=x_1dy_1+\ldots+x_ndy_n=c(x_1ry_1+\ldots+x_nry_n)=0.$$
This is a contradiction.

The second assertion follows from the first assertion: let
$e^2=e\in R$. Then $e\cdot e+(1-e)\cdot (1-e)=1$, and for any
$x\in R$, we have $e(ex-xe)e+(1-e)(ex-xe)(1-e)=0$, whence
$ex-xe=0$.~\hfill$\square$

\textbf{Lemma 2.4.}\label{subgroups}
Let $G$ be a group, $F$ be a field of characteristic $p>0$, and
let $q$ be a prime integer which is not equal to $p$. If the ring
$FG$ is centrally essential, then every $q$-subgroup in $G$ is a
normal commutative subgroup.

\textbf{Proof.} First, let $H$ be a finite $q$-subgroup of the group $G$. Then $|H|=n=q^k$ is a nonzero element of the field $F$ and the element $e_H=\frac 1 n \Sigma_H$ is an idempotent of the ring $FG$. By Lemma 2.3, $e_H$ is a central idempotent. Consequently,
$ge_Hg^{-1}=\frac 1 n \sum_{h\in H}ghg^{-1}=\frac 1 n \sum_{h\in
H}h$ for any $g\in G$. By comparing the coefficients in the both
parts of the last relation, we see that $ghg^{-1}\in H$, i.e., the
subgroup $H$ is normal.

Let $F_0$ be a prime subfield of the field $F$. We consider the
finite ring $F_0H$. By Maschke's theorem, it is isomorphic to some
finite direct product of matrix rings over division rings; in
addition, any finite division ring is a field by Wedderburn's
theorem. We assume that the group $H$ is not commutative. Then one
of the summands of the ring $F_0H$ is the matrix ring of order
$k>1$ over some field; this is impossible, since such a matrix
ring contains a non-central idempotent.

Now let $H$ be an arbitrary $q$-subgroup in $G$. We take any
element $h\in H$ and an arbitrary element $g\in G$. Since $h$
generates the cyclic $q$-subgroup $H_0=\langle h \rangle$, we have
that $ghg^{-1}\in H_0\subseteq H$ for any $g\in G$, i.e., the
subgroup $H$ is normal.

If $x,y\in H$, then the subgroup $H_1=\langle x,y \rangle$ is
finite by
Lemma 2.1(1) and the following Dicman's lemma:\\
if $x_1,\ldots, x_n$ are elements of finite order in an arbitrary
group $G$ and each of the elements $x_1,\ldots, x_n$ has only
finitely many conjugate elements, then there exists a finite
normal subgroup $N$ of the group $G$ containing $x_1,\ldots,x_n$
(see
\cite{Dicman} or \cite[Lemma C, Appendixes]{Lambek}).

By the above, we have $xy=yx$.~\hfill$\square$

In the case of finite groups, we have a more strong assertion
which reduces the study centrally essential group algebras of a
finite group to the study centrally essential group algebras of
finite $p$-groups.

\textbf{Proposition 2.5.}\label{fingr}
Let $|G|=n<\infty$ and $F$ a field of characteristic $p>0$. Then
the following conditions are equivalent.

\textbf{1)} $FG$ is a centrally essential ring.

\textbf{2)} $G=P\times H$, where $P$ is the unique Sylow $p$-subgroup of the group $G$, the group $H$ is commutative, and the ring $FP$ is centrally essential.

\textbf{Proof.}
Let the ring $FG$ be centrally essential. By Lemma 2.4, every
Sylow $q$-subgroup for $q\neq p$ is normal in $G$ and it is
commutative; consequently, the product $H$ of all such subgroups
is a commutative normal subgroup. Let $m=|H|$. We note that
$(m,p)=1$; therefore $m$ is an invertible element of the field
$F$.

We prove that the Sylow $p$-subgroup $P$ of the group $G$ is
normal in $G$.

We consider the following linear mapping $f\colon R\rightarrow R$:
$$f(r)=\frac 1 m \sum_{h\in H}hrh^{-1}.$$
It is obvious that $f(1)=1$ and $f(yry^{-1})=f(r)$ for any $y\in
H$, since the left part and the right part of the relation contain
the same summands. Now we assume that $xy\neq yx$ for some $x\in
P$ and $y\in H$. We set $r=x-yxy^{-1}$. It is directly verified
that $r\neq 0$, but $f(r)=f(x)-f(yxy^{-1})=0$; this contradicts to
Lemma 2.3. Thus, the elements $P$ and $H$ commute, $G=PH$ and
$P\cap H=\{1\}$; consequently, $G=P\times H$. By considering $FG$
as the group ring $(FP)H$, we obtain from Lemma 2.1(1) that $FP$
is a centrally essential ring.

The converse assertion directly follows from Lemma 2.1(2) and the
isomorphism $FG\cong (FP)H$.~\hfill$\square$

\textbf{Proposition 2.6.}\label{nilp2}
Let $G$ be a finite $p$-group and $F$ a field of characteristic
$p$. If $NC(G)\le 2$, then the group ring $FG$ is centrally
essential.

\textbf{Proof.}
We recall that for any subgroup $N$ of the group $G$,  we denote
by $\omega H$ the right ideal of the ring $FG$ generated by the
set $\{1-h|h\in H\}$, we also recall that this right ideal is a
two-sided ideal if and only if the subgroup $H$ is normal. It is
well known (e.g., see cite[Lemma 3.1.6]{Passman}) that the ideal
$\omega G$ is nilpotent in the considered case.

Let $0 \neq x \in FG$. We consider all possible products $x(1-z)$,
where $z\in Z$. If at least one of them (say, $x_1=x(1-z_1)$) is
non-zero, we consider products $x_1(1-z)$ and so on. This process
terminates at some step, i.e., there exists an integer $k\geq 0$
such that $x_k\neq 0$, but $x_k\omega Z=0$ (we assume that
$x_0=x$). Then $x_k\in FG\Sigma_Z$ (see \cite[Lemma
3.1.2]{Passman}). We note that $FG\Sigma_Z\subseteq C(FG)$.
Indeed, if $g,h\in G$, then
$$[g,h\Sigma_Z]=[g,h]\Sigma_Z=gh(1-h^{-1}g^{-1}hg)\Sigma_Z=0,$$
since $h^{-1}g^{-1}hg\in G'\subseteq Z$. Thus, by setting
$c=(1-z_1)\ldots(1-z_k)$ (or $c=1$ for $k=0$), we obtain $c\in
C(FG)$ and $xc=x_k\in C(FG)\setminus\{0\}$, which is
required.~\hfill$\square$

\textbf{Lemma 2.7.}\label{conj}
Let $F$ be a field of characteristic $p$ and let $G$ be a finite $p$-group satisfying the following condition $(*)$:\\
for any non-central element $g\in G$ there exists a
non-trivial subgroup $H$ of the center of the group $G$ with $gH\subseteq g^G$.\\
If $NC(G)>2$, then the ring $R=FG$ is not centrally essential.

\textbf{Proof.}
Let $Z=Z_1(G)$ be the center of the group $G$,~ $K=g^G$ be the
class of conjugate elements of the group $G$,  $|K|>1$,  and let
$H$ be the subgroup from $(*)$.  Then $K=\bigcup_{x\in K}Hx$,
since $a^{-1}Hga=Ha^{-1}ga$ for any $a\in G$. Consequently,
$\Sigma_K=\sum_{i=1}^t \Sigma_{Hx_i}$ for some $x_1,\ldots,x_t\in
K$. Now we note that $(\Sigma_Z)h=\Sigma_Z$ for any $h\in H$;
therefore, $\Sigma_Z\Sigma_H=|H|\Sigma_Z=0$. This implies that
$\Sigma_Z\Sigma_K=0$.

Further, if $NC(G)>2$, then there exists an element $g\in
G\setminus Z_2(G)$. This means that there exists an element $a\in
G$ with $(g,a)\not\in Z$.  We consider the element
$x=g\Sigma_Z\neq 0$. We have
$$
[a,x]=[a,g\Sigma_Z]=(ag-ga)\Sigma_Z=ag(1-(g,a))\Sigma_Z\neq 0,
$$
since $1-(g,a)\not\in \omega Z$. Consequently, $x\not\in C=C(R)$.
An arbitrary element $c\in C$ can be represented as $c=c_0+c_1$,
where $c_0\in FZ$, $c_1=
\sum_{i=0}^s\alpha_i\Sigma_{K_i}$, $K_1,\ldots,K_s$ are the classes
of conjugate elements of the group $G$, and $|K_i|>1$ for all
$i=1,\ldots,s$.  By the above, we have  $xc_1=0$, i.e., $xc=xc_0$.
Since $(\Sigma_Z)z=\Sigma_Z$ for any $z\in Z$, we obtain that
$xc_0=\alpha x$ for some $\alpha\in F$. If $\alpha\neq 0$ and
$xc\in C$, then $x\in C$; this is a contradiction.~\hfill$\square$

\textbf{Lemma 2.8.}\label{centralizer}
If the centralizer of the subgroup $Z_2(G)$ in $G$ is contained in
$Z_2(G)$, then the group $G$ satisfies the condition $(*)$ from
Lemma 2.7.

\textbf{Proof.}
Let $K=g^G$ and $|K|>1$, i.e., $g\not\in Z_1(G)$. If $g\in
Z_2(G)$, then for any element $a\in G$, we have $(g,a)\in Z_1(G)$;
in addition, $(g,a)\neq 1$ for some element $a\in G$. We set
$z=(g,a)$. Then $ga=agz$ and $a^{-1}ga=gz\in K$. Therefore,
$g\langle{z}\rangle\subseteq K$.

Now let $g\not \in Z_2(G)$. Then there exists an element $a\in
Z_2(G)$ such that $z=(g,a)\neq 1$. But $z\in Z_1(G)$ and $ga=agz$,
whence we obtain $g\langle{z}\rangle\subseteq K$.~\hfill$\square$

\textbf{Proposition 2.9.} If $F$ is a field of characteristic $p>0$, then there exists a group $G$ of order $p^5$ such that the group algebra $FG$ is not a centrally essential ring.

\textbf{Proof.} We construct groups which satisfy the conditions of Lemma 2.8.  We consider the cases $p=2$ and $p\neq 2$ separately.

Let $p=2$. We consider the direct product $N$ of the quaternion
group $Q_8=\{\pm 1,\pm i,\pm k\}$ and the cyclic group $\langle a
\rangle$ of order 2 with generator $a$ and the automorphism
$\alpha$ of the group $N$ defined on generators by the relations
$\alpha(i)=j$, $\alpha(j)=i$, $\alpha(a)=(-1)a$. We set
$\Gamma=\langle\alpha\rangle$. We have the semidirect product
$G=N\leftthreetimes\Gamma$ whose elements are considered as
products $x\gamma$, where $x\in N$ and $\gamma \in
\langle\alpha\rangle$, and the operation is defined by the
relation $x\gamma x'\gamma'=x\gamma(x')\gamma\gamma'$. The
elements of the form $x1$ are naturally identified with elements
$x\in N$ and the elements of the form $1\gamma$ are identified
with the elements $\gamma\in \Gamma$. It is directly verified that
$Z_1(G)=\langle-1\rangle$, $Z_2(G)=\langle k,
a\rangle=C_G(Z_2(G))$.

Now we assume that $p>2$. We consider the semidirect product
$N=A\leftthreetimes G$ of the elementary Abelian group  $A$ of
order $p^3$ with generators $a,b,c$ and the group
$\Gamma=\langle\gamma\rangle$, where $\gamma$ is the automorphism
of the group $A$ defined on generators by the relations
$$
\gamma(a)=a,\ \gamma(b)=b,\ \gamma(c)=bc.
$$
It is directly verified that $|N|=p^4$ and any element of the
group $N$ can be uniquely represented in the form of the product
$a^kb^lc^m\gamma^r$, where {$k,l,m,r\in\{0,\ldots,p-1\}$.} We
prove that the mapping $\beta:\{a,b,c,\gamma\}\rightarrow N$,
defined by the relations
$$\beta(a)=a,\ \beta(b)=b,\ \beta(c)=ac,\
\beta(\gamma)=abc\gamma,$$
can be extended to an automorphism $\hat\beta$ of the group $N$.
Indeed, for any $k,l,m,r\in\ZZ$, we set
$$\hat\beta(a^kb^lc^m\gamma^r)=a^kb^la^mc^ma^rb^r(c\gamma)^r=a^{k+m+r}b^{l+\frac {r (r+1)} 2}c^{m+r}\gamma^r.$$
This definition is correct, since $p|\frac{r(r+1)}2$ if $p|r$. It
is directly verified that for any
$k,l,m,r,k',l',m',r'\in\{0,\ldots,p-1\}$ the relations
$$a^kb^lc^m\gamma^r\cdot
a^{k'}b^{l'}c^{m'}\gamma^{r'}=a^{k+k'}b^{l+l'+rm'}c^{m+m'}\gamma^{r+r'}$$
hold. Consequently,
\begin{multline*}\hat\beta(a^kb^lc^m\gamma^r\cdot
a^{k'}b^{l'}c^{m'}\gamma^{r'})=\\
a^{k+k'+m+m'+r+r'}b^{l+l'+rm'+\frac{(r+r')(r+r'+1)}
2}c^{m+m'+r+r'}\gamma^{r+r'}.\end{multline*} On the other hand,
\begin{multline*}
\hat\beta(a^kb^lc^m\gamma^r)\cdot
\hat\beta(a^{k'}b^{l'}c^{m'}\gamma^{r'})=\\
(a^{k+m+r}b^{l+\frac {r(r+1)}2}c^{m+r}\gamma^r)
\cdot(a^{k'+m'+r'}b^{l'+\frac {r'(r'+1)}2}
c^{m'+r'}\gamma^{r'})=\\
 a^{k+m+r+k'+m'+r'}b^{l+\frac {r(r+1)}2+l'+\frac
 {r'(r'+1)}2+r(m'+r')}c^{m+r+m'+r'}\gamma^{r+r'}.
\end{multline*}
It remains to note that we have the following identity
\begin{multline*}
l+\frac {r(r+1)}2+l'+\frac
 {r'(r'+1)}2+r(m'+r')=\\l+l'+rm'+rr'+\frac{r^2+r+{r'}^2+r'} 2=
 l+l'+rm'+\frac{r^2+2rr'+r'^2+r+r'}2=\\l+l'+rm'+\frac{(r+r')(r+r'+1)} 2
\end{multline*}

Now we set $G=N\leftthreetimes \langle\beta\rangle$. It is
directly verified that $Z_1(G)=\langle a,b \rangle$ and
$Z_2(G)=\langle a,b,c
\rangle=C_G(Z_2(G))$.~\hfill$\square$

\textbf{Remark 2.10.} Theorem 1.1 follows from Proposition 2.5,
Proposition 2.6 and Proposition 2.9.

\begin{center}
\textbf{3. Additional Remarks}
\end{center}

\textbf{Remark 3.1.} In connection to
Theorem 1.1(3), we note that a group ring of a finite $p$-group of
nilpotence class 3 can be centrally essential and may also not be
centrally essential. More precisely, we used the computer algebra
system GAP \cite{GAP} to verify that for any group of order
16 which is of nilpotence class
3, its group algebra over the field $\GF(2)$ is centrally
essential.

\textbf{Remark 3.2.} The following semiprimeness criterion of a group ring is well known: the ring $AG$ is semiprime if and only if the ring $A$ is semiprime and the orders of finite normal subgroups of the group $G$ are not zero-divisors in $A$; e.g., see Proposition 8 in  \cite[Appendix]{Lambek}).

\textbf{Proposition 3.3.}\cite[Theorem 2(2), Proposition 4.1]{MT} Any centrally essential semiprime ring $R$ is commutative.

\textbf{Proof.} We verify that for any $r\in R$, the ideal $$r^{-1}C=\{c\in C:\;rc\in C\}$$ is dense in $C$. Indeed, let $d\in C$ and
$dr^{-1}C=0$. If  $dr=0$, then $d\in r^{-1}C$ and $d^2=0$, whence
$d=0$. Otherwise, since $R$ is a centrally essential ring, there
exists an element $z\in C$ with $zdr\in C\setminus\{0\}$. Then
$zd\in r^{-1}C$  and $(zd)^2=0$, whence $zd=0$; this contradicts
to the choice of $z$. The condition ``$r(r^{-1}C)\neq 0$ for any
$r\in R\setminus\{0\}$'' is equivalent to the property that the
ring $R$ is centrally essential. Therefore, $R$ is a right ring of
quotients of the ring $C$ in the sense of \cite[\S 4.3]{Lambek};
consequently, $R$ can be embedded in the complete ring of
quotients of the ring $C$ which is commutative.~\hfill$\square$

\textbf{Proposition 3.4.} Let $A$ be a semiprime ring such that its additive group is torsion-free and let $G$ be an arbitrary group. The ring $AG$ is a centrally essential if and only if the ring $A$ and the group $G$ are commutative.

Proposition 3.4 follows from Remark 3.2 and Proposition 3.3.

\textbf{Remark 3.5.} A. Yu. Ol'shanskii informed the authors of another series of groups satisfying the conditions of Lemma 2.8. Namely, let $p$ be a prime integer and let $G$ be the free 3-generated group of the variety defined by the identities $x^p=1$ and $(x_1,x_2,x_3,x_4)=1$. Then $G/G'$ is an elementary Abelian $p$-group; therefore $G'$ is the Frattini subgroup of the group $G$. If $g\not \in G'$, then can be included $gG'$ in system of free generators of the group $G/G'$; consequently, $g$ can be included in a system consisting of three generators of the group $G$. Since the group $G$ is finite, this generator system is free. Therefore, if $g\in C_G(G')$, then $G$ satisfy the identity $(x_1,x_2,x_3)=1$; this is impossible, since the group $G$ can be mapped onto the group of upper unitriangular matrices of order 4 which does not satisfy this identity. Therefore, $Z_2(G)\supseteq G'\supseteq C_G(G')\supseteq C_G(Z_2(G))$.

\end{document}